# The False Dilemma: Bayesian *vs.* Frequentist[*1]

# 1 Jordi Vallverdú, Ph.D.

Philosophy Dept.

Universitat Autònoma de Barcelona

E-08193 Bellaterra (BCN)

Catalonia - Spain

Jordi.vallverdu@uab.es

**Abstract:** There are two main opposing schools of statistical reasoning, Frequentist and Bayesian approaches. Until recent days, the frequentist or classical approach has dominated the scientific research, but Bayesianism has reappeared with a strong impulse that is starting to change the situation. Recently the controversy about the primacy of one of the two approaches seems to be unfinished at a philosophical level, but scientific practices are giving an increasingly important position to the Bayesian approach. This paper eludes philosophical debate to focus on the pragmatic point of view of scientists' day-to-day practices, in which Bayesian methodology is very useful. Several facts and operational values are described as the core-set for understanding the change.

**Introduction. From dice to propensities.**

When I was developing my PhD research trying to design a comprehensive model to understand scientific controversies and their closures, I was fascinated by statistical problems present in them. The perfect realm of numbers was not able to help to establish clearly causal relationships, nor to agree with one unified statistical conception. Two main schools, Bayesian and frequentist, were fighting each other to demonstrate their own superiority and their own truthful approach. It is an interesting dilemma to analyze. Although I decided then to focus on more general epistemic questions I was

---

[1*] This paper is the result of the communication that I presented on XIIth International Congress of Logic, Methodology and Philosophy of Science, held in Oviedo, $12^{th}$ August 2003.

convinced about the necessity to return to the problem in the future with a thorough study.

Causality has been one of the main topics of the history of philosophical and scientific thought, perhaps the main one: Where do we come from? What is the cause of all things? What will happen? Cosmogonical answers where the first attempt to explain in a causal way the existence of things and beings. The Greek creation myth involved a game of dice between Zeus, Poseidon, and Hades. Later, and beyond myths, Aristotle was the strongest defender of the causal and empirical approach to reality (*Physics*, II, 4-6) although he considered the possibility of chance, especially the problem of the dice game (*On Heavens*, II, 292a30) and probabilities implied in it. But this had nothing to do with the ideas about atomistic chance of Leucippus and Democritus nor Lucrecius' controversial *clinamen's* theory. Hald (1988) affirms the existence of probabilistic thought, not mathematical, in Classical Antiquity. We can also find traces of it in medieval Talmudic and Rabbinical texts, and we know that in 960, the bishop Wibolf of Cambrai calculated 56 diverse ways of playing with 3 dice. *De Vetula*, a Latin poem from XIIIth century, tells us of 216 possibilities. The funny origins of statistical thought: religious man reasoning about games.

In 1494 Luca Paccioli defined the basic principles of algebra and multiplication tables up to 60×60 in his book *Summa de arithmetica, geometria, proportioni e proportionalita*. He posed the first serious statistical problem of two men playing a game called 'balla', which is to end when one of them has won six rounds. However, when they stop playing A has only won five rounds and B three. How should they divide the wager? It would be another 200 years before this problem was solved.

In 1545 Girolamo Cardano wrote the books *Ars magna* (the great art) and *Liber de ludo aleae* (book on games of chance). This was the first attempt to use mathematics to describe statistics and probability, and accurately described the probabilities of throwing various numbers with dice. Galileo expanded on this by calculating probabilities using two dice.

Blaise Pascal (1660) refined the theories of statistics and, with Pierre de Fermat, solved the 'balla' problem of Paccioli. These all paved the way for modern statistics, which essentially began with the use of actuarial tables to determine insurance for merchant ships. In 1662, Antoine Arnauld and Pierre Nicole, publish the influential *La logique ou l'art de penser,* where we can find statistical probabilities. Games and their statistical roots worried people like Cardano, Pascal, Fermat or Huygens (Weatherdord, 1982), although all of them were immersed in a strict mechanistic paradigm. Huygens is considered the first scientist interested in scientific probability, and in 1657 he published *De ratiotiniis in aleae ludo*. Later, De Moivre wrote the influential *De mensura sortis* (1711), and seventy eight years later, Laplace published his *Philosophical Assay About Probability*. In the 1730s, Daniel Bernoulli Jacob's nephew) developed the idea of utility as the mathematical combination of the quantity and perception of risk.

In 1763 an influential paper written by Reverend Thomas Bayes was published posthumously. Richard Price, who was a friend of his, worked on the results of his efforts to find the solution to the problem of computing a distribution for the parameter of a binomial distribution: *An Essay towards solving a Problem in the Doctrine of Chances*. Proposition 9 in the essay represented the main result of Bayes. Degrees of belief are therein considered as a basis for statistical practice. This classical version of Bayesianism had a long history beginning with Bayes and continuing through Laplace to Jeffreys, Keynes and Carnap in the twentieth century. Later, in the 1930's, a new type of Bayesianism appeared, the 'subjective Bayesianism' of Ramsey and De Finetti[2].

At the end of XIXth century, a lot of things were changing in the scientific and philosophical arena. The end of the idea of 'causality' and the conflicts about observation lied at the heart of the debate. Gödel attacked Hilbert's axiomatic approach to mathematics and Bertrand Russell, as clever as ever, told us: "The law of causality (...) is a relic of a bygone age, surviving, like the monarchy, only because it is erroneously supposed to do no harm (...) The principle 'same cause, same effect', which philosophers imagine to be vital to science, is therefore utterly otiose"[3].

In the 1920's arose from the works of Fischer (1922), Neyman and Pearson (1928) the *classic* statistical paradigm: frequentism. They use the relative frequency concept, that is, you must perform one experiment lots of times and measure the proportion where you get a positive result. This proportion,, if you perform the experiment enough times, is the probability. If Neyman and Pearson wrote their first joint paper and presented their approach as *one among alternatives*, Fisher, with his null hypothesis testing[4] gave a different

---

[2] Ramsey (1931), de Finetti (1937) and Savage (1954).

[3] In 1913. Quoted from Suppes (1970): 5.

[4] Nevertheless, this is a controversial concept. See Anderson *et al* (2000). Note that the authors try to find alternatives to null hypothesis testing inside frequentist approach, considering Bayesian methods "computationally difficult and there may continue to be

message: his statistics was the formal solution of the problem of inductive inference (Gigerenzer, 1990: 228).

Philosophers of science like Karl R. Popper were talking about *A World of Propensities.* Nancy Cartwright defends today a probabilistic theory of causation[5]. We live in a world with plenty of uncertainty and risk, because we think that this is the true nature of things. Heisenberg's indetermination principle and the society of risk are different ways by which we understand reality and, at the same time, react to it. Statistics is at the root of our thoughts. But are our thoughts by nature frequentist or Bayesian? Which of the two schools is the better one for scientific practice? This article tries not to address this to philosophers of science but to scientists from diverse scientific fields, from High Energy Physics to Medicine. Beyond philosophers' words, scientific activity makes its own rules. Normative philosophy of science should turn into a prescriptive discipline, in addition to a clarifying and descriptive activity.

These are some of the philosophical and historical[6] aspects of causality but, what about statistical theory and, most importantly, its practice?

**What does 'Bayesian' or 'Frequentist' mean?**

I have omitted in the previous chapter any reference to Bayesian or frequentist approaches. General theory of causality is not necessarily proper to statistical practice, although we must recognize the existence of theories about statistical practice. In fact, two of them are the leading ways to understand several uses of statistics: Bayesian and frequentist approaches.

**Bayesian approach.** This perspective on probabilities, says that a probability is a measure of a person's degree of belief in an event, given the information available. Thus, probabilities refer to a state of knowledge held by an individual, rather than to the properties of a sequence of events. The use of subjective probability in calculations of the expected value of actions is called *subjective expected utility*. There has been a renewal of interest for

---

objections of a fundamental nature to the use of Bayesian methods in strength-of-evidence-assessments and conclusion-oriented, empirical science", p. 921.

[5] Cartwright, N. (1979) 'Causal Laws and Effective Strategies', *Noûs*, 13: 419-437.

[6] For a good and curious history of statistics see Stigler, Stephen M. (1999) *Statistics on the Table: The History of Statistical Concepts and Methods*, USA: Harvard University Press (and his previous works); Salsburg, David (2001) *The lady Tasting Tea: How Statistics Revolutionized Science in the Twentieth Century*, San Francisco: W.H. Freeman and Company; and, finally, two good histories: Hald, A. (1990). *A History of Probability and Statistics and Their Applications before 1750*, NY: Wiley;Hald, A. (1998). *A History of Mathematical Statistics from 1750 to 1930*, NY: Wiley.

Bayesianism since 1954, when L.J. Savage wrote Foundations *of Statistics*. There are a large number of types of Bayesians[7], depending on their attitude towards subjectivity in postulating priors. Recent Bayesian books: Earman (1992), Howson & Urbach (1989), Bernardo & Smith (1996).

**Frequentist approach.** They understand probability as a long-run frequency of a 'repeatable' event and developed a notion of confidence intervals. Probability would be a measurable frequency of events determined from repeated experiments. Reichenbach, Giere or Mayo have defended that approach from a philosophical point of view, referred to by Mayo (1997) as the 'error statistical view" (as opposed to the Bayesian or "evidential-relation view").

**Scientific Activities and the Conductist Approach.**

One of the recurrent arguments against/in favor of one of the two positions (frequentist or Bayesian) consists in saying that a true scientist is always/never frequentist/Bayesian (you can choose between the two possibilities)[8]. It seems to be an epistemological law about statistical practices: "A true scientist never belongs to the opposite statistical school". What I can say is that this is a usual metatheoretical thought about frequentist/Bayesian approaches and their ontological fitting with reality, which is not useful for clarifying the closure of scientific controversies, because they depend on another *kind of values.* We cannot know what happens exactly in scientists' minds, but we can know how they act and, therefore, infer from their actions how they think. Obviously, we suppose and accept cognitive activity for scientists. The question is: at what cost can we introduce cognitive arguments inside the statistical techniques embedded in scientific practices? And when we talk about 'cognition' me must include not only rational aspects of cognition, but also irrational ones[9].

As an initial and simple exercise I tried to search on the Internet for quantitative information about both schools. I typed 'Bayesian' or 'frequentist' on several websites' search engines and found these results: on *Google*, appeared 869.000 results for 'Bayesian' and 18.300 for 'frequentist'; on *Nature's* journal 109 and 2 results, respectively; 82 and 6 on *Science*; finally, on *Medline* 3877 and 124.

---

[7] Ironically, Good, I., J., (1971) told about "46656 kinds of Bayesians", *Amer. Statist.*, 25: 62-63.

[8] As an example, see the complete ideas of Giere, Ronald (1988) *Understanding Scientific Reasoning*, USA: The University of Chicago, p.189 "Are Scientists Bayesian Agents? (...) The overwhelming conclusion is that humans are not bayesian agents", and of B. Efron (1986), "Why isn't everyone a Bayesian?" *American Statistician.* 40: 1-5 or R.D. Cousins (1995), "Why isn't every physicist a Bayesian?" Am. J. Phys. 63: 398. The last two do not need to be quoted.

[9] Paul Thagard (1992) made special software, Explanatory Coherence By Harmony Optimization (ECHO) to analyze scientific knowledge as a set of hypotheses and the relationships established among them:

What does it mean? Is the Bayesian approach the best one? Is it widely used inside hypercommunities but not in the 'common' world? There is a fact: everyday, more and more researchers are using Bayesian methodologies to do their research. Meanwhile, the dominating position of frequentism is weaker and weaker, although still dominating.

But we must go to the core set of the question: what are the new values, beyond scientific values from Merton[10], implied in the choice between the two schools?

**3.1. Formational values:** we will use the words of Bland & Altman (1998): 1160, to illustrate these kinds of values: "Most statisticians have become Bayesians or Frequentists as a result of their choice of university[11]. They did not know that Bayesians and Frequentists existed until it was too late and the choice had been made. There have been subsequent conversions. Some who were taught the Bayesian way discovered that when they had huge quantities of medical data to analyze the frequentist approach was much

---

$$\text{System Coherence} = \sum p_{ij} * a_i * a_j$$

'p' refers to weight, 'a' to acceptation and 'i/j' are the acquired values by the values of a system. Scientists are exposed to hot cognitive variables (such as emotional or motivational variables). Motivation affects the explanatory coherence of scientists. Freedman (1992):332, admits that social processes must be included in the computational models. There is a 'social principle' in the ECHO model: some evidence $E$, produced by a rival in a specific theoretical field appears to have a reduced evidential value. Cognition appeals to something more than mental representations of classic rational values, and must consider the existence of an 'emotional coherence' (as the expanded model HOTCO, for 'hot coherence', Thagard, 2000: 173). See also Thagard (1988, 1999).

[10] The "values" of academic science, as Robert K. Merton wrote in 1942 are: communalism, universality, disinterestedness, originality and skepticism. Merton, R. K. (1942) "The normative structure of science", in The Sociology of Science, N.W. Storer (ed.): 267-278.

[11] The same idea is repeated in a different way by a High Energy physicist, D'Agostini (1998:1): "The intuitive reasoning of physicists in conditions of uncertainty is closer to the Bayesian approach than to the frequentist ideas taught at university and which are considered the reference framework for handling statistical problems". One thing is the theory taught at university, and another one is the true scientific practice.

quicker and more practical, although they remain Bayesian at heart. Some Frequentists have had Damascus road conversions[12] to the Bayesian view. Many practicing statisticians, however, are fairly ignorant of the methods used by the rival camp and too busy to have time to find out". As the epidemiologist Berger says (2003): "practicing epidemiologists are given little guidance in choosing between these approaches apart from the ideological adherence of mentors, colleagues and editors". Giles (2002), talking about members of the Intergovernmental Panel on Climatic Change (IPCC), says that those researchers were suspicious of Bayesian statistics because "these attitudes also stem from the authors' backgrounds", p. 477.

So, the arguments go beyond the ethereal philosophical arena and closer to practical ones. Better opportunities to find a good job is an important argument, and the value of a Bayesian academic training is now accepted: "where once graduate students doing Bayesian dissertations were advised to try not to look too Bayesian when they went on the job market, now great numbers of graduate students try to include some Bayesian flavor in their dissertations to increase their marketability", Wilson (2003): 372.

**3.2. Metaphysical values:** by their writings, we can extract some information about scientist's thoughts. Knowledge is framed by feelings, emotions, facts and, even, faiths. How to consider, then, classical and continuous disputes among the full range of possible positions between realists and subjectivists?

All scientists believe for different reasons, that the constituents of the world have certain dispositions that can be discovered under certain investigative conditions. As expressed by Hacking (1972): 133: "Euler at once retorted that this advice is metaphysical, not mathematical. Quite so! The choice of primitive concepts for inference *is* a matter of 'metaphysics'. The orthodox statistician has made one metaphysical choice and the Bayesian another".

**3.3. Philosophical values (from the scientists' point of view):** to be honest, we must accept that most scientists are not interested in the philosophy of scientific results. But when it an extraordinary fact happens, like a controversy or a paradigm change, they accept it and turn to philosophical ideas (if they have been clearly formulated). Sternberg (2001), writing about the controversies about how to evaluate a diagnostic test, says: "Bayesian methods (...) are all scientifically sound approaches for the evaluation of diagnostic tests in the absence of a perfect gold standard[13], whereas any version of discrepant analysis is not", p. 826. Following the same way of reasoning, the National Academy of Sciences (1993) admitted: "Full Bayesian analyses are often complicated and time-consuming. Moreover, because the data necessary to estimate the component prior probabilities and likelihood ratios may be unavailable, quantitative expression of the assessor's uncertainty is often highly subjective, even if based on expert opinion (...) despite the committee's attempts at objectivity, the interpretation of scientific evidence always retains at least some subjective elements. Use of such 'objective' standards as *P*values, confidence intervals, and relative risks may convey a false sense that such judgments are entirely objective. However judgments about potential sources of bias, although based on sound scientific principles, cannot usually be quantified. This is true even for the scientific 'gold standard' in evaluating causal relationships: the randomized clinical trial" pp.25, 31. In that case, in spite of all that was said, the NAS adopted the Bayesian approach for their analysis. Lilford & Braunholtz (1996:604) also argue: "when the situation is less clear cut (...) conventional statistics may drive decision makers into a corner and produce sudden, large changes in prescribing. The problem does not lie with any of the individual decision makers, but with the very philosophical basis of scientific inference. We propose that conventional statistics should not be used in such cases and that the Bayesian approach is both epistemologically and practically superior". And Spiegelhalter (1999): "There are strong philosophical reasons for using a Bayesian approach".

---

[12] See the curious arguments from a former frequentist: Harrell, Frank E. Jr (2000) *Practical Bayesian data Analysis from a Former Frequentist*, downloadable PDF document at http://hesweb1.ed.virginia.edu/biostat/teaching/bayes.short.course.pdf.

**3.4. Simplicity[14]and cheapness: computerizing statistical thought.** One of the arguments against Bayesian methods says that the Bayesian approach is too complex to apply in day-to-day research. And simplicity is one of the best values for scientific activity[15]. But during the past few

---

[13] Black & Craig (2002): 2653 define 'gold standard' as: "a diagnostic test with 100 per cent sensitivity and specificity". They admit that, frequently, this occurs because of the prohibitive cost or non-existence of a gold standard test. In this situation, rather than using a single imperfect test, multiple imperfect tests may be used to gain an improved prevalence estimate. In general, the results of these tests are correlated, given a subjects' disease status. Bayesian methods are a better solution for these cases than frequentist ones.

[14] See the philosophical reflections about Bayesianism and simplicity of Escoto (2003).

[15] Gigerenzer, Gerd (1989) "We Need Statistical Thinking, Not Statistical Rituals", *Behavioral and Brain Sciences*, 21, 2 (1998): 199-200 is a

years a large amount of Bayesian software programs have appeared which have changed the situation: now it is easy, fast and cheap to implement the Bayesian approach in experimental practices. Programs like BACC, [B/D], BOA, BUGS (Bayesian inference using Gibbs sampling, and WinBUGS), MINITAB, EPIDAT, FIRST BAYES, HYDRA, STATA, SAS, S-Plus and others, some of them available as freeware, make possible an efficient use of Bayesian methods in several scientific fields. Their flexibility helps to incorporate multiple sources of data and of uncertainty within a single coherent composite model. Until the 1980's, the potential for the application of Bayesian methods was limited by the technical demands placed on the investigator. Over the past fifteen years these limitations have been substantially reduced by innovations in scientific computing (faster computer processors)[16] and drastic drops in the cost of computing (Editorial *BMJ*, 1996). These changes and an increase in the number of statisticians trained in Bayesian methodology are encouraging the new status of Bayesianism (Tan, 2001).

Medicine is, perhaps, the scientific field in which Bayesian analysis is being more intensively applied (Szolovits, 1995; Grunkemeir & Payne, 2002). Two trends, evidence-based medicine and Bayesian statistics are changing the practice of contemporary medicine. As Ashby & Smith (2000) tells us: "Typically the analysis from such observational studies [those of epidemiology] is complex, largely because of the number of covariates. Probably for this reason, Bayesian applications in epidemiology had to wait for the recent explosion in computer power, but are now appearing in growing numbers", p. 3299[17]. We must also take into account the expert (Bayesian) systems; a topic developed on section 4.2.

The development of Markov chain Monte Carlo (MCMC) computation algorithms, now permit fitting models with incredible realistic complexity[18].

---

very critical of the complexities of frequentist methods. Downloadable at http://www.mpib-berlin.mpg.de/dok/full/gg/ggwnsbabs/ggwnsbabs.html. A similarly aggressive work: Rindskopf, D. (1998) "Null-hypothesis tests are not completely stupid, but Bayesian statistics are better", *Behavioral and Brain Sciences*, 21: 215-216.

[16] See NAS (1991).

[17] See also Breslow (1990) and Ahsby & Hutton (1996).

[18] When we study models for multiple comparisons, we can see that frequentists adjust Multiple Comparison Procedures (MCP) considering intersection of multiple null hypotheses. They also advocate for a control of the familywise error-rate (FWE). So, "Bayesians will come closer to a frequentist per-comparison or to a FEW approach depending on the credibility they attach to the family of (null) hypotheses being tested (...) the Bayesian is closer to the per-comparisonist", Berry (1999): 216.

The Bayesian approach has received a great impulse from MCMC models (Dunson, 2001; Carli & Louis, 2000; Gelman *et al* 1996). MCMC procedures are also extremely flexible and constitute the primary factor responsible for the increased use and visibility of Bayesian methods in recent years.

**3.5. Ethical values:** we can find an appeal to ethical values as parts of arguments about both schools. Wilson (2003) affirms that Bayesian methods are a more ethical approach to clinical trials and other problems. On the contrary, Fisher (1996) affirms that "Ethical difficulties may arise because of the differing types of belief", especially during Randomized Clinical Trials (the Phase III Trials in the FDA model).

From the history of standard literature on ethics in medical research, man can infer the great value of prior beliefs in clinical trials. And the key concept is 'uncertainty': "Subjective opinions are typically not included in the background material in a clinical trial protocol, but as they are often a driving force behind the existence of a protocol, and as uncertainty is deemed to be ethically important, documentation will be useful. Without documentation it may be difficult to determine whether uncertainty exists. (...) There are compelling ethical reasons that uncertainty should be present before a clinical trial is undertaken" (Chaloner & Rhame, 2001: 591 and 596). When uncertainty is introduced in the reasoning procedures, the quantification of prior beliefs and, therefore, the use of Bayesian methodologies, seems to be an operationally and ethically better decision.

**3.6. Better fitting for results:** Berger (2003), proposes using both models and studying case by case their possibilities: "based on the philosophical foundations of the approaches, Bayesian models are best suited to addressing hypotheses, conjectures, or public-policy goals, while the frequentist approach is best suited to those epidemiological studies which can be considered 'experiments', i.e. testing constructed sets of data". Usually, we find no such equitable position.

But this is not a theoretical question but a practical one: Bayesian methods work better than frequentist. Therefore, Bayesian methods are increasing their application range, although it does not always mean that there are more 'true Bayesians'. As Wilson (2003) explains:" their methodological successes [from Bayesian] have indeed impressed many within the field and without, but those who have adopted the Bayesian methods have often done so with-out adopting the Bayesian philosophy"[19]. As Popper or Lakatos[20] could say: "Bayesian methods solve problems better than frequentist ones". And practical success usually means the theory's success. Look to the history of science: Copernicus astronomical tables were better than those of Ptolomeus and if at first, were accepted as an instrument, in a later they were considered as a true representation of reality.

So, The Scientific Information and Computing Center at CIBA-GEIGY's Swiss headquarters in Basle moved towards the systematic use of Bayesian methods not so much as a result of theoretical conviction derived from philosophical debates, but rather as a pragmatic response to the often experienced inadequacy of traditional approaches to deal with the problems with which CIBA-GEIGY statisticians were routinely confronted (Racine *et al*, 1986). An example: clinical trials made by pharmaceutical industries are usually Bayesian (Estey & Thall, 2003) although such methods are not easily implemented (Wang *et al*, 2002).

Bayesian methods are ideally suited to dealing with multiples sources of uncertainty, and risk assessment must include a lot of them: one experiment can be affected by several terms like sex, age, occupation, skill of technician, number of specimens, time of sampling, genetic background, source of intake... So, according to an epidemiologist, Dunson (1991): 1225: "Bayesian approaches to the analysis of epidemiological data represent a powerful tool for interpretation of study results[21] and evaluation of hypotheses about exposure-disease relations. These tools allow one to consider a much broader class of conceptual and mathematical models than would have been possible using non-Bayesian approaches"[22]. Grunkmeier & Payne (2002: 1901), talking about surgery enumerate several advantages of Bayesian statistics applied to it: "(1) providing direct probability statements – which are what most people wrongly assume they are getting from conventional statistics; (2) formally incorporating previous information in statistical inference of a data set, a natural approach which follows everyday reasoning; and (3) flexible, adaptive research designs allowing multiple examination of accumulating study data". The Bayesian approach is more efficient at unifying and calculating multilevel causal relationships[23].

---

[19] See the Editorial from *British Medical Journal* (1996), "most people find Bayesian probability much more akin to their own thought processes (...) The areas in which there is most resistance to Bayesian methods are those were the frequentist paradigm took root in the 1940s to 1960s, namely clinical trials and epidemiology. Resistance is less strong in areas where formal inference is not so important, for example during phase I and II trials, which are concerned mainly with safety and dose finding".

[20] Popper, Karl R. (1963). Conjectures and Refutations: The Growth of Scientific Knowledge. London: Routledge and Kegan Paul. Popper, Karl R. (1972). Objective Knowledge: An Evolutionary Approach. Oxford: Oxford University Press. Popper, Karl R. (1959). The Logic of Scientific Discovery. New York: Basic Books, Inc. Lakatos, I. 'Falsification and the Methodology of Scientific Research Programmes', in Lakatos, I & Musgrove, A. (eds). *Criticism and the Growth of Knowledge*. Cambridge University Press, Cambridge, 1970; Lakatos, I. *The Methodology of Scientific Research Programmes*, (ed. J. Worrall & G. Currie). Cambridge University Press, 1978; Lakatos, I & Musgrove, A. (eds). *Criticism and the Growth of Knowledge*. Cambridge University Press, Cambridge, 1970.

[21] For experimental Bayesian design see Chaloner & Verdinelli (1995).

**3.7. Diffusion of science: guidelines.** At the core of science remains information communication. By the process of writing and communicating his/her results, a scientist is at the same time evaluated (through *peer review*) and judged (by his/her colleagues). All the norms implied in the guidelines, define a trend in 'good' scientific practices[24]. And those groups who control the communication channels can make sure that special kinds of ideas are never allowed. Therefore, design and control of communication channels is something crucial for the interest of a community.

The frequentist approach has dominated statistics journals all through XXth Century but, recently, Bayesians are gaining more and more power. As Wilson (2003): 372, says: "Bayesians have successfully and extensively published in JASA and other prominent journals, bringing their methods into the spotlight where they cannot be ignored". It is not only a question of general perception but also of radical changes in the bases of the epistemic frame. The International Committee of Medical Journal Editors, wrote the *Uniform Requirements for Manuscripts Submitted to Biomedical Journals*, which you can consult at http://www.icmje.org, where they specified for sta-tistical norms: "Avoid relying solely on statistical hypothesis testing, such as the use of P values, which fail to convey important quantitative information"[25].

---

[22] Freedman (1996) remarks that epidemiological studies cannot be usually made without external information. He also affirms that the choice of a $p-$value, like P¡0.05, implies the inclusion of a subjective factor in the evaluation of the experimental results.

[23] Thagard (1999) has offered a very powerful conceptual framework to understand scientific explanations of diseases with his idea of "causal network instantiation" (p. 114). According to him: "causal networks are not simple schemas that are used to provide single causes for effects, but they instead describe complex mechanisms of multiple interacting factors", p. 115-116. But Thagard is no Bayesian: he pursues another line of explanation which he considers better suited to psychological reasoning: *explanatory coherence.* (Ibid. p. 65-66).

[24] We must also include publication bias such as the quicker publishing of papers of studies with positive results than those with null or negative findings, (Dickersin *et al*, 2002).

Spiegelhalter (1999) reflects that: "Current international guidelines for statistical submissions to drug regulatory authorities state that 'the use of Bayesian and other approaches may be considered when the reasons for their use are clear and when the resulting conclusions are sufficiently robust'"[26].

So, these new trends 'accepted' as the new axiological frame for statistical research have changed the weight of both schools: while frequentist models are decreasing their expansion, Bayesian ones are being employed in an increasing number of situations. Basañez (2004) has explained the reasons for this gradual shift: practical, theoretical and philosophical.

**Cognition and statistics:**

**Looking into the scientist's (human) mind.** If we have looked to external activity in the previous chapter, now we must analyze the internal or cognitive activity of human beings. The first and more important sense for humans is the visual capacity. The latest studies about human visual processes (Geisler & Kersten, 2002) show that Bayesian explanations fit better than frequentist when we must explain how we process visual information and react properly to it. The Bayesian approach seems to be optimal to explain in the broader biological context, plasticity, learning and natural selection[27]. Perception and inference work in a Bayesian way (Knill *et al*, 1996). If we consider the statistical properties of natural environments and how these interact in the process of natural selection to determine the design of perceptual and cognitive systems, we must accept that the Bayesian framework captures and generalizes, in a formal way, many of the important ideas of other approaches to perception and cognition (Geisler & Diehl, 2003). Computational (Marr, 1982) and evolutionary (Pinker, 1997) studies are well explained and unified by the Bayesian framework (Liu & Kersten, 1998; Shiffrin & Steyvers, 1997; Legge, Klitz & Tjan, 1997; Dosher & Lu, 2000; Stankiewicz, Legge & Mansfield, 2000). But human evolution can be employed to demonstrate the opposite arguments: "in his evaluation of evidence, man is apparently not a conservative Bayesian: he is not a Bayesian at all" (Kahneman & Tversky, 1982).

---

[25] We must also recognize that the use of statistical methodologies in medical research is highly controversial, beyond the Bayesian-frequentist dilemma (Altman *et al,* 2002).

[26] See http://www.findarticles.com/p/articles/mi_m0999/is_7208_319/ai_55721117 [electronic resource], accessed August, $1^{st}$ 2004.

[27] Nevertheless, we can find completely opposite opinions. As an example, Fisher (1996) says: "Humans do not and cannot behave in a Bayesian manner", p. 424. And he also justifies his ideas appealing to human evolution!

Neural networks[28] are difficult to study because of their complexity and in the Bayesian approach these issues can be handled in a natural and consistent way. Several problems which appeared in the standard neural network methods can be solved by the Bayesian approach. For example, "the unknown degree of complexity and the resulting model is handled by defining vague (non-informative) priors for the hyper-parameters that determine the model complexity, and the resulting model is averaged over all model complexities weighted by their posterior probability giving the data sample. (...) The Bayesian analysis yields posterior predictive distributions for any variables of interest, making the computation of confidence intervals possible", Lampinen & Vehtari (2001):7.

An important moment in the controversies about both schools was the book written by Howson and Urbach in 1989 (see also the paper of 1991) about hypothesis evaluation and inferences using the Bayesian approach. They affirm that the Bayesian view starts off acknowledging that subjective assessment of likelihood is an important part of theory selection and construction, and makes it part of the philosophy of science. The power of scientific reasoning then, results not from some elusive objective logic of discovery, but because our innate inference abilities lead from observation of evidence to beliefs that follow probability calculus, and hence our sense of increasing credibility tends to reflect greater likelihood of a theory making accurate predictions. It follows that, our beliefs can be measured as probabilities, and probabilities can be used to confirm theories. In that case, novel observations should have and do have special importance in theory construction. The authors introduce probability calculus in simple algebraic terms and discuss its application to the philosophy of science.

**Robot Minds: Expert Systems and the AI.** One of the Artificial Intelligence trends is the design of expert systems. These are intended to replicate the decision making of a human expert within narrowly defined domains[29]. Such domains are highly specialized: Logic Theorist was the first

---

[28] And the "Bayesian Networks", introduced in the late 1980's and early 1990's, representing an important advance in probabilistic reasoning for artificial intelligence, as you have seen, and expert systems. See Glymour (2003).

[29] And such systems have three primary components: a knowledge base, decision rules, and an inference engine (Kurzweil, 1990: 292).

expert system and was made in 1955 by Herbert Simon and Allan Newell as a logic theorist. This expert system discovered most of the principles written by Russell & Whitehead in *Principia Mathematica*. We also have Heuristic Dendral (mass spectrography, 1967), Macsyma (indefinite integrals solver, 1967), Macsyma (mathematics teacher, 1967), Internist (internal medicine, 1970), Mycin (blood infections specialist, 1974), Prospector (geology, 1978), and so on until the present day[30]. The engineering AI approach to computational philosophy of science is allied with "¡¡android epistemology¿¿, the epistemology of machines that may or may not be built like humans" (Thagard, in Bynum 1998: 52). Most current expert systems apply Bayesian probability to their studies. One of them has discovered previously unsuspected fine structure in the infrared spectra of stars (Cheeseman, 1990)

But what is the connection between AI's expert systems and statistics and why is it mostly developed in a Bayesian way? The keystone of all this is the idea of *uncertainty*[31]. We humans, must decide a lot of actions without complete data about facts we are analyzing, that is, we must take decisions in the light of uncertain knowledge about a situation. In 1947, von Neumann and Morgenstern developed an axiomatic framework called *utility theory*, founded on the pre-existing axioms of probability theory (Horvitz, 1993). Utility theory provides a formal definition of preference and of rational decisions under uncertainty, and its axioms define a measure of preference called *utility*. Then, utility theory affirms that people should make decisions that have optimal average, or expected, utility. Those expected values fit with the idea of personal or *subjective probability*. The use of subjective probability in calculations of the expected value of actions is called *subjective expected utility* (SEU).

In the mid-1940s, 'operation analysis' was developed, matured into the modern discipline of *operations research* (OR), and used after World War II to analyze a broad spectrum of civilian problems.

SEU and OR were developed closely and led to the emergence, in the early 1950s, of management science and decisions analysis. The 1960s represented the creation of the first expert systems, the 1970s were the realm of vision research and, finally, the 1980s represented the maturation of rules-based expert systems. Today, we must realize that: "subjective probability methodology has proved extremely successful for gathering the probabilities needed by an expert system" (Horvitz, 1993:31).

---

[30] In 1975 Dendral discovered a new rule to identify organic molecules, which no human mind had ever thought. Prospector discovered in 1982 a big molybdenum deposit.

[31] Understood, basically, as a situation with a lack of information. For more accurate definitions of 'uncertainty' see Zimmermann (2000).

Machines are certainly not humans, but they work better and in a more similar way to humans when they use Bayesian principles (Horvitz *et al*, 1988). Good results with these expert systems (such as AI learning to recognize objects, Viola, 1996) are giving powerful reasons in favor of Bayesianism. As an example, the Stanford researcher Paté-Cornell, wrote a paper (2002) on fusion intelligence and the Bayesian approach, trying to apply that method to the USA Intelligence Services, and avoid another attack after the attack suffered on September $11^{th}$, 2001. Bayesianism has reached new and more influential levels of application and is gaining ground.

**Framing values: conclusions about theories and uses.**

My old Webster's Dictionary has its own definition of "dilemma": "**1.** a situation requiring a choice between equally undesirable alternatives. **2.** any difficult or perplexing situation or problem. **3.** *Logic.* a form of syllogism in which the major premise is formed of two or more hypothetical propositions and the minor premise is an exhaustive disjunctive proposition, as ¡¡If A, then B; if C then D. Either A or C. Therefore, either B or D¿¿".

It seems clear that we have not been talking about logic relationships inside statistical controversies. Therefore, the third definition is not of interest. The second one seems to be closer to the aims of this paper: the analysis of a complex problem for which there is no obvious solution. Finally, the last shall be first, the first definition is the core of this paper: Does the Bayesian *vs.* frequentist dilemma constitute a difficult choice 'between equally undesirable alternatives'? Are we forced to die for our rational criteria like Buridan's donkey?

At a metatheoretical level, that is philosophy, the debate is still open and more and more complex. But that is not the level of analysis we have considered as crucial for the solution of the debate. We talk about scientific practices in which are involved both statistical approaches. And when scientists work, they take decisions continuously.

We have shown a new range of values that constitute part of the statistical axiology. These are non-epistemic values, but shape the underlying framework of research epistemology. Academic training, ease of use, powerful infrastructures, cognitive fitting, ethics, metaphysical options, cheapness, and better results, are the arguments to decide in favor of either one of the two approaches. Perhaps these are not the values which theoreticians would have chosen, but are the real values which appear when we look at scientists' practices and reflections.

We don't know if the prediction made by Bruno de Finetti, that it would take until the year 2020 for the Bayesian view of statistics to completely prevail will be accurate. This is another question, far from our interests and methodology. I have indicated several values that makes it possible to choose between both approaches.

A clear fact is that Bayesian analysis is widely used in a variety of fields, from the pioneering field of medicine to engineering, image processing, expert systems, decision analysis, psychological diagnoses (Meehl & Rosen, 1955), criminal investigations (Sullivan & Delaney, 1982), for presenting evidence in court (Feinberg & Schervish, 1986; Matthew, 1994; Mossman, 2000), gene sequencing, financial predictions, neural networks or epidemiological studies. If we return to the classic paper of Winkler (1974) "Why are experimental psychologists (and others) reluctant to use Bayesian inferential procedures in practice?"[32], we will read: "this state of affairs appears to be due to a combination of factors including philosophical conviction, tradition, statistical training, lack of 'availability', computational difficulties, reporting difficulties, and perceived resistance by journal editors". Well, all these factors (non-epistemic values) are now not against but in favor of the Bayesian approach.

Is the solution to unify as a synthesis both approaches (Berger *et al*, 1997), like a synthesis solution to a dualistic problem? Could a hybrid method of inference satisfy both camps? Is the Likelihood approach a third alternative? (Senn, 2003). But this is, once more, a philosophical question.

Finally, we must consider the existence of a really fundamental question: how to make decisions based on evidence. And we find a basic problem: there are several decision levels with their own individual exigencies regarding what is considered as evidence. If we talk about decision making in health controversies, we should consider several levels like: decision making for patients (diagnosis), decision making for individual patients (interventions), decision making about studies (start from prior beliefs and data monitoring), decision making for pharmaceutical companies and public policy decision making (Ashby & Smith, 2000). But these multi-criteria analyses can be found in other scientific fields, such as forestry (Kangas & Kangas, 2004). And we find another set of problems present in both approaches when they are applied to controversial scientific practices, such as those of risk assessment: difficulties in establishing clear relationships, the significant sample, data interpretation, cognitive paradoxes (Simpson, Ellsberg, St. Petersbourg, Base-Rate Fallacy,...), the idea of evidence at multiple levels[33], use of both models,...

---

[32] Quoted by Lecoutre, Bruno at http://www.stat.auckland.ac.nz/~iase/publications/5/leco0735.pdf [electronic document]. Accessed in July $30^{th}$, 2004.

Considering the previous arguments, we must admit that the dilemma, understood as a choice between equally undesirable alternatives, is a false dilemma. We have enough judgment elements to decide rationally[34] between one of two approaches, and so do scientists from diverse fields, whose words we have reproduced here. To understand these decisions better we have enumerated a new set of values that needs to be included in a richer and sounder scientific axiology.

---

[33] In the herbicide 2,4,5-T controversy we found: "Petitioners demand sole reliance on *scientific* facts, on evidence that reputable scientific techniques certify as certain. Typically, a scientist will not so certify evidence, unless the probability of error, by standard statistical measurement, is less than 5%. That is, scientific fact is at least 95% certain. (...) Agencies are not limited to scientific fact, to 95% certainties. Rather, they have at least the same fact-finding powers as a jury, particularly when, as here, they are engaged in rule-making" (Jasanoff, 1994: 51). But, at the same time: "Typically a scientist will not...certify evidence unless the probability of error, by standard statistical measurement, is less than 5%. That is, scientific fact is at least 95% certain. Such certainty has never characterized the judicial or the administrative process...the standard of ordinary civil litigation, a preponderance of evidence, demands only 51% certainty. A jury may weigh conflicting evidence and certify as adjudicative (although not scientific) fact that which it believes is more likely than not...Inherently, such a standard is flexible; inherently, it allows the fact-finder to assess risks, to measure probabilities, to make subjective judgments. Nonetheless, the ultimate finding will be treated, at low as fact...The standard before administrative agencies is no less flexible. Agencies are not limited to scientific fact, to 95% certainties...we must deal with the terminology of law, not science", Miller (1980): 75-76 (Miller made extracts from trial "Ethyl Corp. V. EPA, 541F .2d 1 (1976), p.28,.58").

[34] Rationality must be understood as a complex activity well modeled by researchers like Kuhn. See the interesting work of Salmon, Wesley C. (1990) Rationality *and Objectivity in Science or Tom Kuhn Meets Tom Bayes*, in Savage, C. (ed.) *Minnesota Studies in the Philosophy of Science*, Vol. XIV, USA: University of Minnesota Press 175-204.